\documentclass[12pt,noamsfonts]{amsart}

\usepackage{amsfonts}
\usepackage{amssymb}
\newcommand{\Spec}{\operatorname{Spec}}
\newcommand{\Proj}{\operatorname{Proj}}
\newcommand{\codim}{\operatorname{codim}}

\newtheorem{question}{Question}

\theoremstyle{remark}
\newtheorem*{solution}{Solution}

\begin{document}

\title[M2 and Schemes]{Macaulay~2 and the Geometry of Schemes} 
\author[Smith]{Gregory G.~Smith}
\address{Department of Mathematics, University of California,
Berkeley, California, 94720} 
\email{{\tt ggsmith@math.berkeley.edu}, {\tt bernd@math.berkeley.edu}}
\author[Sturmfels]{Bernd Sturmfels}

\maketitle

This tutorial illustrates how to use Dan Grayson and Mike Stillman's
computer algebra system, {\sc Macaulay2}, to study schemes.  The
examples are taken from the homework for an algebraic geometry class
given at the University of California, Berkeley in the fall of 1999.
This graduate course was taught by the second author with assistance
from the first author.  Our choice of problems follows the material in
David Eisenbud and Joe Harris' textbook {\em The Geometry of Schemes}.
In fact, four of our ten problems are taken directly from their
exercises.

%%----------------------------------------------------------
%%  QUESTION #1
%%----------------------------------------------------------
\subsection*{Distinguished open sets}

We begin with a simple example involving the definition of an affine
scheme; see section~I.1.4 in Eisenbud and Harris (1999).  This example
also indicates some of the subtleties involved in working with
arithmetic schemes.

\begin{question} 
Let $R = {\mathbb Z}[x,y,z]$ and $X = \Spec(R)$; in other words, $X$
is affine $3$-space over the integers.  Let $f = x$ and consider the
basic open subset $X_{f} \subset X$.
\begin{enumerate}
\item[(a)] If $e_{1} = x+y+z$, $e_{2} = xy+xz+yz$ and $e_{3} = xyz$
are the elementary symmetric functions then the set $\left\{ X_{e_{i}}
\right\}_{1 \leq i \leq 3}$ is an open cover of $X_{f}$.
\item[(b)] If $p_{1} = x+y+z$, $p_{2} = x^{2}+y^{2}+z^{2}$ and $p_{3}
= x^{3}+y^{3}+z^{3}$ are the power sum symmetric functions then
$\left\{ X_{p_{i}} \right\}_{1 \leq i \leq 3}$ is {\em NOT} an open
cover of $X_{f}$.
\end{enumerate}
\end{question}

\begin{solution}
By Lemma~I-6 in Eisenbud and Harris (1999), it suffices to show that
$e_{1}$, $e_{2}$ and $e_{3}$ generate the unit ideal in $R_{f}$.  This
is equivalent to showing that $x^{m}$ belongs to the $R$-ideal
$\langle e_{1}, e_{2}, e_{3} \rangle$ for some $m \in {\mathbb N}$.
In particular, the saturation $\big( \langle e_{1}, e_{2}, e_{3}
\rangle : x^{\infty} \big)$ is the unit ideal if and only if $\left\{
X_{e_{i}} \right\}_{1 \leq i \leq 3}$ is an open cover of $X_{f}$.
{\sc Macaulay2} allows us to work with homogenous ideals over
${\mathbb Z}$ and we obtain: {\scriptsize
\begin{verbatim}
    i1 : R = ZZ[x, y, z];
\end{verbatim}
\begin{verbatim}
    i2 : elementaryBasis = ideal(x+y+z, x*y+x*z+y*z, x*y*z);
\end{verbatim}
\begin{verbatim}
    o2 : Ideal of R
\end{verbatim}
\begin{verbatim}
    i3 : saturate(elementaryBasis, x)
\end{verbatim}
\begin{verbatim}
    o3 = ideal 1
\end{verbatim}
\begin{verbatim}
    o3 : Ideal of R
\end{verbatim}}
Similarly, to prove that $\left\{ X_{p_{i}} \right\}_{1 \leq i \leq
3}$ is not an open cover of $X_{f}$, it is enough to show that $\big(
\langle p_{1}, p_{2}, p_{3} \rangle : x^{\infty} \big)$ is not the
unit ideal.  We compute this saturation:
{\scriptsize
\begin{verbatim}
    i4 : powerSumBasis = ideal(x+y+z, x^2+y^2+z^2, x^3+y^3+z^3);
\end{verbatim}
\begin{verbatim}
    o4 : Ideal of R
\end{verbatim}
\begin{verbatim}
    i5 : saturate(powerSumBasis, x)
\end{verbatim}
\begin{verbatim}
                                2           2
    o5 = ideal (6, x + y + z, 2y  - y*z + 2z , 3y*z)
\end{verbatim}
\begin{verbatim}
    o5 : Ideal of R
\end{verbatim}}
\noindent However, working over the field ${\mathbb Q}$, we find that
$\big( \langle p_{1}, p_{2}, p_{3} \rangle : x^{\infty} \big)$ is the
unit ideal.  {\scriptsize
\begin{verbatim}
    i6 : S = QQ[x, y, z];
\end{verbatim}
\begin{verbatim}
    i7 : powerSumBasis = ideal(x+y+z, x^2+y^2+z^2, x^3+y^3+z^3);
\end{verbatim}
\begin{verbatim}
    o7 : Ideal of S
\end{verbatim}
\begin{verbatim}
    i8 : saturate(powerSumBasis, x)
\end{verbatim}
\begin{verbatim}
    o8 = ideal 1
\end{verbatim}
\begin{verbatim}
    o8 : Ideal of S
\end{verbatim}}
\end{solution}

%%----------------------------------------------------------
%%  QUESTION #2
%%----------------------------------------------------------
\subsection*{Irreducibility}

The study of complex semisimple Lie algebras gives rise to an
important family of algebraic varieties called nilpotent orbits.  To
illustrate one of the properties appearing in section I.2.1 of
Eisenbud and Harris (1999), we examine the irreducibility of a
particular nilpotent orbit.

\begin{question} 
Let $X$ be the set of complex $3 \times 3$ matrices which are
nilpotent.  Show that $X$ is an irreducible algebraic variety.
\end{question}

\begin{solution}
A $3 \times 3$ matrix $M$ is nilpotent if and only if its minimal
polynomial divides ${\mathsf T}^{k}$, for some $k \in {\mathbb N}$.
Since each irreducible factor of the characteristic polynomial of $M$
is also a factor of the minimal polynomial, we conclude that the
characteristic polynomial of $M$ is ${\mathsf T}^{3}$.  It follows
that the coefficients of the characteristic polynomial (except for the
leading coefficient which is $1$) of a generic $3 \times 3$ matrix
define the algebraic variety $X$.

To show $X$ is irreducible over $\mathbb{C}$, it is enough to
construct a rational parameterization of $X$; see Proposition~4.5.6 in
Cox, Little, and O'Shea (1996).  To achieve this, observe that
$\mathrm{GL}_{n}(\mathbb{C})$ acts on $X$ by conjugation.  Jordan's
canonical form theorem implies that there are exactly three orbits;
one for each of the following matrices:
\[
N_{0} = \left[ \begin{smallmatrix} 0 & 0 & 0 \\ 0 & 0 & 0 \\ 0 & 0 &
0 \end{smallmatrix} \right] \, , \,
N_{1} = \left[ \begin{smallmatrix} 0 & 1 & 0 \\ 0 & 0 & 0 \\ 0 & 0 &
0 \end{smallmatrix} \right]
\text{ and } 
N_{2} = \left[ \begin{smallmatrix} 0 & 1 & 0 \\ 0 & 0 & 1 \\ 0 & 0 &
0 \end{smallmatrix} \right] \, .
\]
Each orbit is defined by a rational parameterization, so it suffices
to show that the closure of the orbit containing $N_{2}$ is the entire
variety $X$.  In {\sc Macaulay2}, this calculation can be done as
follows: 
{\scriptsize
\begin{verbatim}
    i1 : S = QQ[t, y_0..y_8, a..i, MonomialOrder => Eliminate 10];
\end{verbatim}
\begin{verbatim}    
    i2 : N2 = (matrix {{0,1,0},{0,0,1},{0,0,0}}) ** S
\end{verbatim}
\begin{verbatim}
    o2 = {0} | 0 1 0 |
         {0} | 0 0 1 |
         {0} | 0 0 0 |
\end{verbatim}
\begin{verbatim}
                 3       3
    o2 : Matrix S  <--- S
\end{verbatim}
\begin{verbatim}    
    i3 : G = genericMatrix(S, y_0, 3, 3);
\end{verbatim}
\begin{verbatim}
                 3       3
    o3 : Matrix S  <--- S
\end{verbatim}}
\noindent To determine the entries in $\, G N_{2} G^{-1}$, we use the
classical adjoint to construct the inverse of the generic matrix $G$.
{\scriptsize
\begin{verbatim}
    i4 : adj = (G,i,j) -> (
              n := degree target G;
              m := degree source G;
              (-1)^(i+j)*det(submatrix(G, {0..(i-1), (i+1)..(n-1)}, 
                        {0..(j-1), (j+1)..(m-1)}))
              );
\end{verbatim}
\begin{verbatim}    
    i5 : classicalAdjoint = (G) -> (
              n := degree target G;
              matrix table(n, n, (i, j) -> adj(G, j, i))
              );
\end{verbatim}
\begin{verbatim}
    i6 : numerators = G*N2*classicalAdjoint(G);
\end{verbatim}
\begin{verbatim}
                 3       3
    o6 : Matrix S  <--- S
\end{verbatim}
\begin{verbatim}    
    i7 : D = det(G);
\end{verbatim}
\begin{verbatim}
    i8 : M = genericMatrix(S, a, 3, 3);
\end{verbatim}
\begin{verbatim}
                 3       3
    o8 : Matrix S  <--- S
\end{verbatim}}
\noindent The entries in $G N_{2} G^{-1}$ give a rational
parameterization of the orbit generated by $N_{2}$.  We give an
``implicit representation'' of this variety by using elimination
theory; see section 3.3 in Cox, Little, and O'Shea (1996).  
{\scriptsize
\begin{verbatim}
    i9 : elimIdeal = minors(1, (D*id_(S^3))*M-numerators) + ideal(1 - D*t);
\end{verbatim}
\begin{verbatim}
    o9 : Ideal of S
\end{verbatim}
\begin{verbatim}    
    i10 : closureOfOrbit = ideal(selectInSubring(1, gens gb elimIdeal));
\end{verbatim}
\begin{verbatim}
    o10 : Ideal of S
\end{verbatim}}
\noindent Finally, we check that the closure of this orbit is equal to
$X$ scheme-theoretically.  
{\scriptsize
\begin{verbatim}
    i11 : X = ideal submatrix((coefficients({0}, 
                       det(M - t*id_(S^3))))_1, {1,2,3})
\end{verbatim}
\begin{verbatim}
    o11 = ideal (a + e + i, b*d - a*e + c*g + f*h - a*i - e*i, 
               - c*e*g + b*f*g + c*d*h - a*f*h - b*d*i + a*e*i)
\end{verbatim}
\begin{verbatim}
    o11 : Ideal of S
\end{verbatim}
\begin{verbatim}
    i12 : closureOfOrbit == X
\end{verbatim}
\begin{verbatim}
    o12 = true
\end{verbatim}}
\end{solution}

%%----------------------------------------------------------
%%  QUESTION #3
%%----------------------------------------------------------
\subsection*{Singular Points}

Section~I.2.2 in Eisenbud and Harris (1999) provides the definition of
a singular point of a scheme.  In our third question, we study the
singular locus of a family of elliptic curves.  Section~V.3 in
Eisenbud and Harris (1999) also contains related material.

\begin{question} 
Consider a general form of degree $3$ in ${\mathbb Q}[x,y,z]$;
$$
F = ax^{3} + bx^{2}y + cx^{2}z + dxy^{2} + exyz + fxz^{2} + gy^{3} +
hy^{2}z + iyz^{2} + jz^{3}.
$$
Give necessary and sufficient conditions in terms of $a,
\ldots, j$ for the cubic curve $\Proj\big( {\mathbb Q}[x,y,z] /
\langle F \rangle \big)$ to have a singular point.
\end{question}

\begin{solution}
A time consuming elimination gives the degree $12$ polynomial which
defines the singular locus of a general form of degree $3$.  This can
be done in {\sc Macaulay2} as follows.  We have not displayed the
output {\tt o6}, as this discriminant has $2040$ terms in the $10$
variables $a, \ldots, j$.  
{\scriptsize
\begin{verbatim}
    i1 : S = QQ[x, y, z, a..j, MonomialOrder => Eliminate 2];
\end{verbatim}
\begin{verbatim}
    i2 : F = a*x^3+b*x^2*y+c*x^2*z+d*x*y^2+e*x*y*z+
         f*x*z^2+g*y^3+h*y^2*z+i*y*z^2+j*z^3;
\end{verbatim}
\begin{verbatim}
    i3 : partials = submatrix(jacobian matrix{{F}}, {0..2}, {0})
\end{verbatim}
\begin{verbatim}
    o3 = {1} | 3x2a+2xyb+y2d+2xzc+yze+z2f |
         {1} | x2b+2xyd+3y2g+xze+2yzh+z2i |
         {1} | x2c+xye+y2h+2xzf+2yzi+3z2j |
\end{verbatim}
\begin{verbatim}
                 3       1
    o3 : Matrix S  <--- S
\end{verbatim}
\begin{verbatim}
    i4 : singularities = ideal(partials) + ideal(F);
\end{verbatim}
\begin{verbatim}
    o4 : Ideal of S
\end{verbatim}
\begin{verbatim}
    i5 : elimDiscr = ideal selectInSubring(1, gens gb singularities);
\end{verbatim}
\begin{verbatim}
    o5 : Ideal of S
\end{verbatim}
\begin{verbatim}
    i6 : elimDiscr = substitute(elimDiscr, {z => 1});
\end{verbatim}
\begin{verbatim}
    o6 : Ideal of S
\end{verbatim}}
\noindent There is also a simple and more useful determinantal formula
for this discriminant. It is a specialization of the formula (2.8) in
section~3.2 in Cox, Little, and O'Shea (1998): 
{\scriptsize
\begin{verbatim}
    i7 : hessian = det submatrix(jacobian ideal partials, {0..2}, {0..2});
\end{verbatim}
\begin{verbatim}
    i8 : A = (coefficients({0,1,2}, submatrix(
                       jacobian matrix{{F}}, {0..2}, {0})))_1; 
\end{verbatim}
\begin{verbatim}
                 3       6
    o8 : Matrix S  <--- S
\end{verbatim}
\begin{verbatim}
    i9 : B = (coefficients({0,1,2}, submatrix( 
                       jacobian matrix{{hessian}}, {0..2}, {0})))_1;
\end{verbatim}
\begin{verbatim}
                 3       6
    o9 : Matrix S  <--- S
\end{verbatim}
\begin{verbatim}
    i10 : detDiscr = ideal det (A || B);
\end{verbatim}
\begin{verbatim}
    o10 : Ideal of S
\end{verbatim}
\begin{verbatim}
    i11 : detDiscr == elimDiscr
\end{verbatim}
\begin{verbatim}
    o11 = true
\end{verbatim}}
\end{solution}

%%----------------------------------------------------------
%%  QUESTION #4
%%----------------------------------------------------------
\subsection*{Fields of Definition}

Schemes over non-algebraically closed fields arise in number theory.
Our solution to Exercise~II-6 in Eisenbud and Harris (1999) indicates
one technique for working over a number field in {\sc Macaulay2}.

\begin{question}
An inclusion of fields $K \hookrightarrow L$ induces a map
$\mathbb{A}_{L}^{n} \to \mathbb{A}_{K}^{n}$.  Find the images in
${\mathbb A}_{\mathbb{Q}}^{2}$ of the following points of ${\mathbb
A}_{\overline{\mathbb{Q}}}^{2}$ under this map.
\begin{enumerate}
\item[(a)]  $\langle x - \sqrt{2}, y - \sqrt{2} \rangle$;
\item[(b)]  $\langle x - \sqrt{2}, y - \sqrt{3} \rangle$;  
\item[(c)]  $\langle x - \zeta, y - \zeta^{-1} \rangle$ where $\zeta$ is
a $5$-th root of unity ;
\item[(d)]  $\langle \sqrt{2}x- \sqrt{3}y \rangle$;
\item[(e)]  $\langle \sqrt{2}x- \sqrt{3}y-1 \rangle$.
\end{enumerate}
\end{question}

\begin{solution}
The images can be determined by (1) replacing coefficients not
belonging to $K$ with indeterminates, (2) adding the minimal
polynomials of these coefficients to the given ideal in ${\mathbb
A}_{\overline{\mathbb{Q}}}^{2}$ and (3) eliminating the new
indeterminates. Here are the five examples:
{\scriptsize
\begin{verbatim}
    i1 : S = QQ[a, b, x, y, MonomialOrder => Eliminate 2];
\end{verbatim}
\begin{verbatim}
    i2 : Ia = ideal(x-a, y-a, a^2-2);
\end{verbatim}
\begin{verbatim}
    o2 : Ideal of S
\end{verbatim}
\begin{verbatim}
    i3 : ideal selectInSubring(1, gens gb Ia)
\end{verbatim}
\begin{verbatim}
                        2
    o3 = ideal (x - y, y  - 2)
\end{verbatim}
\begin{verbatim}
    o3 : Ideal of S
\end{verbatim}
\begin{verbatim}
    i4 : Ib = ideal(x-a, y-b, a^2-2, b^2-3);
\end{verbatim}
\begin{verbatim}
    o4 : Ideal of S
\end{verbatim}
\begin{verbatim}
    i5 : ideal selectInSubring(1, gens gb Ib)
\end{verbatim}
\begin{verbatim}
                 2       2
    o5 = ideal (y  - 3, x  - 2)
\end{verbatim}
\begin{verbatim}
    o5 : Ideal of S
\end{verbatim}
\begin{verbatim}
    i6 : Ic = ideal(x-a, y-a^4, a^4+a^3+a^2+a+1);
\end{verbatim}
\begin{verbatim}
    o6 : Ideal of S
\end{verbatim}
\begin{verbatim}
    i7 : ideal selectInSubring(1, gens gb Ic)
\end{verbatim}
\begin{verbatim}
                          2    2               3    2
    o7 = ideal (x*y - 1, x  + y  + x + y + 1, y  + y  + x + y + 1)
\end{verbatim}
\begin{verbatim}
    o7 : Ideal of S
\end{verbatim}
\begin{verbatim}
    i8 : Id = ideal(a*x+b*y, a^2-2, b^2-3);
\end{verbatim}
\begin{verbatim}
    o8 : Ideal of S
\end{verbatim}
\begin{verbatim}
    i9 : ideal selectInSubring(1, gens gb Id)
\end{verbatim}
\begin{verbatim}
                2   3  2
    o9 = ideal(x  - -*y )
                    2
\end{verbatim}
\begin{verbatim}
    o9 : Ideal of S
\end{verbatim}
\begin{verbatim}
    i10 : Ie = ideal(a*x+b*y-1, a^2-2, b^2-3);
\end{verbatim}
\begin{verbatim}
    o10 : Ideal of S
\end{verbatim}
\begin{verbatim}
    i11 : ideal selectInSubring(1, gens gb Ie)
\end{verbatim}
\begin{verbatim}
                 4     2 2   9  4    2   3  2   1
    o11 = ideal(x  - 3x y  + -*y  - x  - -*y  + -)
                             4           2      4
\end{verbatim}
\begin{verbatim}
    o11 : Ideal of S
\end{verbatim}}
\end{solution}

%%----------------------------------------------------------
%%  QUESTION #5
%%----------------------------------------------------------
\subsection*{Multiplicity}

The multiplicity of a zero-dimensional scheme $X$ at a point $p \in X$
is defined to be the length of the local ring $\mathcal{O}_{X,p}$.
Unfortunately, we cannot work directly in the local ring in {\sc
Macaulay2}. What we can do, however, is to compute the multiplicity by
computing the degree of the component of $X$ supported at $p$; see
page~66 in Eisenbud and Harris (1999).

\begin{question}
What is the multiplicity of the origin $(0,0,0)$ as a zero of the
polynomial equations
\[
x^{5}+y^{3}+z^{3} =  x^{3}+y^{5}+z^{3} = x^{3}+y^{3}+z^{5} =  0 \, ?
\]
\end{question}

\begin{solution}
If $I$ is the ideal generated by $x^{5}+y^{3}+z^{3}$,
$x^{3}+y^{5}+z^{3}$ and $x^{3}+y^{3}+z^{5}$ in ${\mathbb Q}[x,y,z]$,
then the multiplicity of the origin is
\[
\dim_{\mathbb{Q}} \frac{\mathbb{Q}[x,y,z]_{\langle x,y,z \rangle}}
{I \cdot \mathbb{Q}[x,y,z]_{\langle x,y,z \rangle}} \, .
\]
It follows that the multiplicity is the vector space dimension of the
ring $\mathbb{Q}[x,y,z] / \varphi^{-1}(I \cdot
\mathbb{Q}[x,y,z]_{\langle x,y,z \rangle})$ where $\varphi \colon
\mathbb{Q}[x,y,z] \to \mathbb{Q}[x,y,z]_{\langle x,y,z \rangle}$ is
the natural map.  Moreover, we can express this using ideal quotients:
\[
\varphi^{-1}(I \cdot \mathbb{Q}[x,y,z]_{\langle x,y,z \rangle}) =
\big(I : (I : \langle x,y,z \rangle^{\infty})\big) \, .
\]
Carrying out this calculation in {\sc Macaulay2}, we obtain:
{\scriptsize
\begin{verbatim}
    i1 : S = QQ[x, y, z];
\end{verbatim}
\begin{verbatim}    
    i2 : I = ideal(x^5+y^3+z^3, x^3+y^5+z^3, x^3+y^3+z^5);
\end{verbatim}
\begin{verbatim}
    o2 : Ideal of S
\end{verbatim}
\begin{verbatim}
    i3 : multiplicity = degree(I : saturate(I))
\end{verbatim}
\begin{verbatim}
    o3 = 27
\end{verbatim}}
\end{solution}

%%----------------------------------------------------------
%%  QUESTION #6
%%----------------------------------------------------------
\subsection*{Flat Families}

Non-reduced schemes arise naturally as the flat limit of a family of
reduced schemes. Exercise~III-68 in Eisenbud and Harris (1999)
illustrates how a family of skew lines in ${\mathbb P}^{3}$ gives a
double line with an embedded point.

\begin{question}
Let $L$ and $M $ be the lines in ${\mathbb P}^{3}_{k[t]}$ given by
$x=y=0$ and $x-tz = y+t^{2}w =0$ respectively.  Show that the flat
limit as $t \to 0$ of the union $L \cup M$ is the double line $x^{2} =
y = 0$ with an embedded point of degree $1$ located at the point
$(0:0:0:1)$.
\end{question}

\begin{solution}
We find the flat limit by saturating the intersection ideal:
{\scriptsize
\begin{verbatim}
    i1 : PP3 = QQ[x, y, z, w];
\end{verbatim}
\begin{verbatim}
    i2 : S = QQ[t, x, y, z, w];
\end{verbatim}
\begin{verbatim}
    i3 : phi = map(PP3, S, 0 | vars PP3 );
\end{verbatim}
\begin{verbatim}
    o3 : RingMap PP3 <--- S
\end{verbatim}
\begin{verbatim}
    i4 : L = ideal(x, y);
\end{verbatim}
\begin{verbatim}
    o4 : Ideal of S
\end{verbatim}
\begin{verbatim}
    i5 : M = ideal(x-t*z, y-t^2*w);
\end{verbatim}
\begin{verbatim}
    o5 : Ideal of S
\end{verbatim}
\begin{verbatim}
    i6 : X = intersect(L, M);
\end{verbatim}
\begin{verbatim}
    o6 : Ideal of S
\end{verbatim}
\begin{verbatim}
    i7 : Xzero = trim phi substitute(saturate(X, t), {t => 0})
\end{verbatim}
\begin{verbatim}
                      2        2
    o7 = ideal (y*z, y , x*y, x )
\end{verbatim}
\begin{verbatim}
    o7 : Ideal of PP3
\end{verbatim}}
\noindent This is the union of a double line and an embedded point of
degree $1$.
{\scriptsize
\begin{verbatim}    
    i8 : use PP3;
\end{verbatim}
\begin{verbatim}
    i9 : intersect(ideal(x^2, y), ideal(x, y^2, z))
\end{verbatim}
\begin{verbatim}
                      2        2
    o9 = ideal (y*z, y , x*y, x )
\end{verbatim}
\begin{verbatim}
    o9 : Ideal of PP3
\end{verbatim}
\begin{verbatim}
    i10 : degree( ideal(x^2, y) / ideal(x, y^2, z))
\end{verbatim}
\begin{verbatim}
    o10 = 1
\end{verbatim}}
\end{solution}

%%----------------------------------------------------------
%%  QUESTION #7
%%----------------------------------------------------------
\subsection*{B\'{e}zout's Theorem}

B\'{e}zout's theorem (Theorem~III-78 in Eisenbud and Harris, 1999)
fails without the Cohen-Macaulay hypothesis.  Following Exercise
III-81 in Eisenbud and Harris (1999), we illustrate this in {\sc
Macaulay2}.

\begin{question}
Find irreducible closed subvarieties $X$ and $Y$ in ${\mathbb P}^{4}$
with
\begin{eqnarray*}
\codim(X \cap Y) & = & \codim(X) +\codim(Y) \text{ and }\\
\deg(X \cap Y) & > & \deg(X) \cdot \deg(Y) \, .
\end{eqnarray*}
\end{question}

\begin{solution}
We show that the assertion holds when $X$ is the cone over the
nonsingular rational quartic curve in ${\mathbb P}^{3}$ and $Y$ is a
two-plane passing through the vertex of the cone.  The computation is
done as follows: 
{\scriptsize
\begin{verbatim}
    i1 : S = QQ[a, b, c, d, e];
\end{verbatim}
\begin{verbatim}
    i2 : quarticCone = trim minors(2, 
              matrix{{a, b^2, b*d, c}, {b, a*c, c^2, d}})
\end{verbatim}
\begin{verbatim}
                            3      2     2    2    3    2
    o2 = ideal (b*c - a*d, c  - b*d , a*c  - b d, b  - a c)
\end{verbatim}
\begin{verbatim}
    o2 : Ideal of S
\end{verbatim}
\begin{verbatim}
    i3 : plane = ideal(a, d);
\end{verbatim}
\begin{verbatim}
    o3 : Ideal of S
\end{verbatim}
\begin{verbatim}
    i4 : codim quarticCone + codim plane == codim (quarticCone + plane)
\end{verbatim}
\begin{verbatim}
    o4 = true
\end{verbatim}
\begin{verbatim}
    i5 : (degree quarticCone) * (degree plane)
\end{verbatim}
\begin{verbatim}
    o5 = 4
\end{verbatim}
\begin{verbatim}
    i6 : degree (quarticCone + plane)
\end{verbatim}
\begin{verbatim}
    o6 = 5
\end{verbatim}}
\end{solution}

%%----------------------------------------------------------
%%  QUESTION #8
%%----------------------------------------------------------
\subsection*{Constructing Blow-ups}

The blow-up of a scheme $X$ along a subscheme $Y$ can be constructed
from the Rees algebra associated to the ideal sheaf of $Y$ in $X$; see
Theorem~IV-22 in Eisenbud and Harris (1999).  Gr\"{o}bner basis
techniques allow one to express the Rees algebra in terms of
generators and relations.  We demonstrate this by solving
Exercise~IV-43 in Eisenbud and Harris (1999).

\begin{question}
Find the blow-up of the affine plane ${\mathbb A}^{2} = \Spec\big(
k[x,y] \big)$ along the subscheme defined by $\langle x^{3}, xy, y^{2}
\rangle$.
\end{question}

\begin{solution}
We first provide a general function which given an ideal and a list of
variables returns the ideal of relations for the Rees algebra.
{\scriptsize
\begin{verbatim}
    i1 : blowUpIdeal = method();
\end{verbatim}
\begin{verbatim}
    i2 : blowUpIdeal(Ideal, List) := (I, L) -> (
             r := numgens I;
             S := ring I;
             kk := coefficientRing S;
             n := numgens S;
             y := symbol y;
             St := kk[t, gens S , y_1..y_r, MonomialOrder => Eliminate 1];
             phi := map(St, S, submatrix(vars St, {1..n}));
             F := phi gens I;
             local J;
             J = ideal apply(1..r, j -> y_j - t*F_(0, (j-1)));
             J = ideal selectInSubring(1, gens gb J);
             if (#L < r) then error "not enough variables";
             R := kk[gens S, L];
             theta := map(R, St, 0 | vars R);
             theta J
         );
\end{verbatim}}
\noindent Applying the function to our specific case yields: 
{\scriptsize
\begin{verbatim}                 
    i3 : S = QQ[x, y];             
\end{verbatim}
\begin{verbatim}
    i4 : I = ideal(x^3, x*y, y^2);
\end{verbatim}
\begin{verbatim}
    o4 : Ideal of S
\end{verbatim}
\begin{verbatim}
    i5 : blowUpIdeal(I, {A, B, C})
\end{verbatim}
\begin{verbatim}
                              2         2          3     2
    o5 = ideal (y*B - x*C, x*B  - A*C, x B - y*A, x C - y A)
\end{verbatim}
\begin{verbatim}
    o5 : Ideal of QQ [x, y, A, B, C]
\end{verbatim}}
\end{solution}

%%----------------------------------------------------------
%%  QUESTION #9
%%----------------------------------------------------------
\subsection*{A Classic Blow-up}

We consider the blow-up of the projective plane ${\mathbb P}^{2}$ at a
point.  Many related examples appear in section IV.2.2 of Eisenbud and
Harris (1999).

\begin{question}
Show that the following varieties are isomorphic.
\begin{enumerate}
\item[(a)] the image of the rational map from ${\mathbb P}^{2}$ to
${\mathbb P}^{4}$ given by
\[
(r:s:t) \mapsto (r^{2}:s^{2}:rs:rt:st) \, ;
\]
\item[(b)] the blow-up of the plane ${\mathbb P}^{2}$ at the point
$(0:0:1)$;
\item[(c)] the determinantal variety defined by the $2 \times 2$
minors of the matrix
\[
\begin{bmatrix} a & c & d \\ b & d & e \end{bmatrix}
\]
where ${\mathbb P}^{4} = \Proj\big( k[a,b,c,d,e] \big)$.
\end{enumerate}
This surface is called the {\em cubic scroll} in ${\mathbb P}^{4}$.
\end{question}

\begin{solution}
We find the ideal in part~(a) by elimination theory.
{\scriptsize
\begin{verbatim}
    i1 : PP4 = QQ[a, b, c, d, e];
\end{verbatim}
\begin{verbatim}
    i2 : S = QQ[r, s, t, A..E, MonomialOrder => Eliminate 3 ];
\end{verbatim}
\begin{verbatim}
    i3 : I = ideal(A - r^2, B - s^2, C - r*s, D - r*t, E - s*t);
\end{verbatim}
\begin{verbatim}
    o3 : Ideal of S
\end{verbatim}
\begin{verbatim}
    i4 : phi = map(PP4, S, (matrix{{0, 0, 0}}**PP4) | vars PP4)
\end{verbatim}
\begin{verbatim}
    o4 = map(PP4,S,{0, 0, 0, a, b, c, d, e})
\end{verbatim}
\begin{verbatim}
    o4 : RingMap PP4 <--- S
\end{verbatim}
\begin{verbatim}
    i5 : surfaceA = phi ideal selectInSubring(1, gens gb I)
\end{verbatim}
\begin{verbatim}
                                             2
    o5 = ideal (c*d - a*e, b*d - c*e, a*b - c )
\end{verbatim}
\begin{verbatim}
    o5 : Ideal of PP4
\end{verbatim}}
\noindent We determine the surface in part~(b) by constructing the
blow-up of ${\mathbb P}^{2}$ in ${\mathbb P}^{2} \times {\mathbb
P}^{1}$ and then projecting its Segre embedding from ${\mathbb P}^{5}$
into ${\mathbb P}^{4}$.  Notice that its image under the Segre map lies
on a hyperplane in ${\mathbb P}^{5}$.  
{\scriptsize
\begin{verbatim}
    i6 : R = QQ[t, x, y, z, u, v, MonomialOrder => Eliminate 1];
\end{verbatim}
\begin{verbatim}
    i7 : blowUpIdeal = ideal selectInSubring(1, gens gb ideal(u-t*x, v-t*y))
\end{verbatim}
\begin{verbatim}
    o7 = ideal(y*u - x*v)
\end{verbatim}
\begin{verbatim}
    o7 : Ideal of R
\end{verbatim}
\begin{verbatim}
    i8 : PP2xPP1 = QQ[x, y, z, u, v];
\end{verbatim}
\begin{verbatim}
    i9 : psi = map(PP2xPP1, R, 0 | vars PP2xPP1);
\end{verbatim}
\begin{verbatim}
    o9 : RingMap PP2xPP1 <--- R
\end{verbatim}
\begin{verbatim}
    i10 : blowUp = PP2xPP1 / psi(blowUpIdeal);
\end{verbatim}
\begin{verbatim}
    i11 : PP5 = QQ[A, B, C, D, E, F];
\end{verbatim}
\begin{verbatim}
    i12 : segre = map(blowUp, PP5, matrix{{x*u, y*u, z*u, x*v, y*v, z*v}});
\end{verbatim}
\begin{verbatim}
    o12 : RingMap blowUp <--- PP5
\end{verbatim}
\begin{verbatim}
    i13 : ker segre
\end{verbatim}
\begin{verbatim}
                                    2
    o13 = ideal (B - D, C*E - D*F, D  - A*E, C*D - A*F)
\end{verbatim}
\begin{verbatim}
    o13 : Ideal of PP5
\end{verbatim}
\begin{verbatim}
    i14 : theta = map( PP4, PP5, matrix{{a, c, d, c, b, e}})
\end{verbatim}
\begin{verbatim}
    o14 = map(PP4,PP5,{a, c, d, c, b, e})
\end{verbatim}
\begin{verbatim}
    o14 : RingMap PP4 <--- PP5
\end{verbatim}
\begin{verbatim}
    i15 : surfaceB = trim theta ker segre
\end{verbatim}
\begin{verbatim}
                                              2
    o15 = ideal (c*d - a*e, b*d - c*e, a*b - c )
\end{verbatim}
\begin{verbatim}
    o15 : Ideal of PP4
\end{verbatim}}
\noindent Finally, we compute the surface in part~(c) and apply a
permutation of the variables to obtain the desired isomorphisms
{\scriptsize
\begin{verbatim}
    i16 : determinantal = minors(2, matrix{{a, c, d},{b, d, e}}) 
\end{verbatim}
\begin{verbatim}
                                              2
    o16 = ideal (- b*c + a*d, - b*d + a*e, - d  + c*e)
\end{verbatim}
\begin{verbatim}
    o16 : Ideal of PP4
\end{verbatim}
\begin{verbatim}
    i17 : sigma = map(PP4, PP4, matrix{{d, e, a, c, b}});
\end{verbatim}
\begin{verbatim}
    o17 : RingMap PP4 <--- PP4
\end{verbatim}
\begin{verbatim}
    i18 : surfaceC = sigma determinantal
\end{verbatim}
\begin{verbatim}
                                              2
    o18 = ideal (c*d - a*e, b*d - c*e, a*b - c )
\end{verbatim}
\begin{verbatim}
    o18 : Ideal of PP4
\end{verbatim}
\begin{verbatim}
    i19 : surfaceA == surfaceB 
\end{verbatim}
\begin{verbatim}
    o19 = true
\end{verbatim}
\begin{verbatim}
    i20 : surfaceB == surfaceC
\end{verbatim}
\begin{verbatim}
    o20 = true
\end{verbatim}}
\end{solution}

%%----------------------------------------------------------
%%  QUESTION #10
%%----------------------------------------------------------
\subsection*{Fano Schemes}

Our last example concerns the family of Fano schemes associated to a
flat family of quadrics.  We solve Exercise~IV-69 in Eisenbud and
Harris (1999).

\begin{question}
Consider the one-parameter family of quadrics 
\[ 
Q = V(tx^{2}+ty^{2}+tz^{2}+w^{2}) \subseteq {\mathbb P}^{3}_{k[t]} =
\Proj\big(k[t][x,y,z,w]\big) \, .
\]
As the fiber $Q_{t}$ tends to the double plane $Q_{0}$, what is the
flat limit of the Fano scheme $F_{1}(Q_{t})$ of lines lying on these
quadric surfaces~?
\end{question}

\begin{solution}
We first make the homogeneous coordinate ring of the ambient
projective $3$-space and the ideal of our family of quadrics. 
{\scriptsize
\begin{verbatim}
    i1 : PP3overBase = QQ[t, x, y, z, w];
\end{verbatim}
\begin{verbatim}
    i2 : Qt = ideal(t*x^2+t*y^2+t*z^2+w^2);
\end{verbatim}
\begin{verbatim}
    o2 : Ideal of PP3overBase
\end{verbatim}}
\noindent We construct an indeterminate line in
$\mathbb{P}_{\mathbb{Q}[t]}^{3}$ by adding parameters $u, v$ and two
points {\footnotesize $(A \colon B \colon C \colon D)$} and
{\footnotesize $(E \colon F \colon G \colon H)$}.  The map $\phi$
sends the variables to the coordinates of the general point on this
line.  
{\scriptsize
\begin{verbatim}
    i3 : S = QQ[t, u, v, A..H];
\end{verbatim}
\begin{verbatim}
    i4 : phi = map(S, PP3overBase, matrix{{t}} | 
              u*matrix{{A, B, C, D}}+v*matrix{{E, F, G, H}});
\end{verbatim}
\begin{verbatim}
    o4 : RingMap S <--- PP3overBase
\end{verbatim}}
\noindent The indeterminate line is contained in our family of
quadrics if and only if $\phi(tx^{2}+ty^{2}+tz^{2}+w^{2})$ vanishes
identically in $u$ and $v$.  Thus, we extract the coefficients of $u$
and $v$.
{\scriptsize
\begin{verbatim}
    i5 : imageOfQt = phi Qt;
\end{verbatim}
\begin{verbatim}
    o5 : Ideal of S
\end{verbatim}
\begin{verbatim}
    i6 : coeff = (coefficients({1,2}, gens imageOfQt))_1;
\end{verbatim}
\begin{verbatim}
                 1       3
    o6 : Matrix S  <--- S
\end{verbatim}}
\noindent We no longer need the variables $u$ and $v$.
{\scriptsize
\begin{verbatim}
    i7 : Sprime = QQ[t, A..H];
\end{verbatim}
\begin{verbatim}
    i8 : coeff = substitute(coeff, Sprime);
\end{verbatim}
\begin{verbatim}
                      1            3
    o8 : Matrix Sprime  <--- Sprime
\end{verbatim}
\begin{verbatim}
    i9 : Sbar = Sprime / (ideal coeff);
\end{verbatim}}
\noindent To move to the Grassmannian over $\mathbb{Q}[t]$, we
introduce a polynomial ring in $6$ new variables corresponding to the
minors of the matrix $\left[ \begin{smallmatrix} A & B & C & D \\ E &
F & G & H \end{smallmatrix} \right]$.  The map $\psi$ sends the new
variables $a, \ldots f$ to the appropriate minor, regarded as elements
of ${\tt Sbar}$.
{\scriptsize
\begin{verbatim}
    i10 : PP5overBase = QQ[t, a..f];
\end{verbatim}
\begin{verbatim}
    i11 : psi = matrix{{Sbar_"t"}} | substitute(
               exteriorPower(2, matrix{{A, B, C, D}, {E, F, G, H}}), Sbar);
\end{verbatim}
\begin{verbatim}
                     1          7
    o11 : Matrix Sbar  <--- Sbar
\end{verbatim}
\begin{verbatim}
    i12 : fanoOfQt = trim ker map(Sbar, PP5overBase, psi);
\end{verbatim}
\begin{verbatim}
    o12 : Ideal of PP5overBase
\end{verbatim}}
\noindent We next determine the limit as $t$ tends to $0$.  
{\scriptsize
\begin{verbatim}    
    i13 : fanoOfQ0 = trim substitute(saturate(
              fanoOfQt, t), {t => 0})
\end{verbatim}
\begin{verbatim}
                                            2   2   2 
    o13 = ideal (e*f, d*f, d*e, a*e + b*f, d , f , e , c*d - b*e + a*f, 
\end{verbatim}
\begin{verbatim}
                                        2    2    2
                 b*d + c*e, a*d - c*f, a  + b  + c )
\end{verbatim}
\begin{verbatim}
    o13 : Ideal of PP5overBase
\end{verbatim}
\begin{verbatim}
    i14 : degree( ideal(d, e, f, a^2+b^2+c^2) / fanoOfQ0)
\end{verbatim}
\begin{verbatim}
    o14 = 2
\end{verbatim}}
\noindent We see that $F_{1}(Q_{0}) $ is supported on the plane conic
$\langle d, e, f, a^{2}+b^{2}+c^{2} \rangle$ and $F_{1}(Q_{0})$ is not
reduced --- it has multiplicity two.

From section~IV.3.2 in Eisenbud and Harris (1999), we know that
$F_{1}(Q_{1})$ is the union of two conics lying in complementary
planes.  We verify this as follows: 
{\scriptsize
\begin{verbatim}
    i15 : fanoOfQ1 = trim substitute(saturate(
              fanoOfQt, t), {t => 1})
\end{verbatim}
\begin{verbatim}
                             2    2    2
    o15 = ideal (a*e + b*f, d  + e  + f , c*d - b*e + a*f, b*d + c*e, 
\end{verbatim}
\begin{verbatim}
                             2    2    2                         2    2 
                 a*d - c*f, c  + e  + f , b*c + d*e, a*c - d*f, b  - e , 
\end{verbatim}
\begin{verbatim}
                             2    2
                 a*b + e*f, a  - f )
\end{verbatim}
\begin{verbatim}
    o15 : Ideal of PP5overBase
\end{verbatim}
\begin{verbatim}
    i16 : fanoOfQ1 == intersect(ideal(c-d, b+e, a-f, d^2+e^2+f^2), 
               ideal(c+d, b-e, a+f, d^2+e^2+f^2))
\end{verbatim}
\begin{verbatim}
    o16 = true
\end{verbatim}}
\noindent Thus, $F_{1}(Q_{0})$ is the double conic obtained when the
two conics in $F_{1}(Q_{1})$ move together.
\end{solution}

%%----------------------------------------------------------
%%  REFERENCES
%%----------------------------------------------------------

\end{document}